\documentclass{amsart}
\usepackage{amssymb}
\usepackage{graphicx}
\usepackage{lineno}

\newtheorem{theorem}{Theorem}

\newcommand{\A}{{\mathcal A}}
\newcommand{\U}{{\mathcal U}}

\newcommand{\IC}{{\mathbb C}}

\newcommand{\ID}{{\mathbb D}}

\def\be{\begin{equation}}
\def\ee{\end{equation}}
\newcommand{\bthm}{\begin{theorem}}
\newcommand{\ethm}{\end{theorem}}
\newcommand{\beqq}{\begin{eqnarray*}}
\newcommand{\eeqq}{\end{eqnarray*}}

\begin{document}

\title[Coefficients of the inverse of functions for the subclass of  $\mathcal U (\lambda)$]{Coefficients of the inverse of functions for the subclass of the class  $\boldsymbol{\mathcal U (\lambda)}$}

\author[M. Obradovi\'{c}]{Milutin Obradovi\'{c}}
\address{Department of Mathematics,
Faculty of Civil Engineering, University of Belgrade,
Bulevar Kralja Aleksandra 73, 11000, Belgrade, Serbia}
\email{obrad@grf.bg.ac.rs}

\author[N. Tuneski]{Nikola Tuneski}
\address{Department of Mathematics and Informatics, Faculty of Mechanical Engineering, Ss. Cyril and Methodius
University in Skopje, Karpo\v{s} II b.b., 1000 Skopje, Republic of North Macedonia.}
\email{nikola.tuneski@mf.edu.mk}



\begin{abstract}
Let $\A$ be the class of functions $f$ that are analytic in the unit disk $\ID$ and normalized such that $f(z)=z+a_2z^2+a_3z^3+\cdots$. Let $0<\lambda\le1$ and
\[ \U(\lambda) = \left\{ f\in\A: \left |\left (\frac{z}{f(z)} \right )^{2}f'(z)-1\right | < \lambda,\, z\in\ID \right\}. \]
In this paper sharp upper bounds of the first three coefficients of the inverse function $f^{-1}$ are given in the case when
$$\frac{f(z)}{z}\prec \frac{1}{(1-z)(1-\lambda z)}.$$
\end{abstract}

\maketitle

Let ${\mathcal A}$ denote the family of all analytic functions
in the unit disk $\ID := \{ z\in \IC:\, |z| < 1 \}$   satisfying the normalization
$f(0)=0= f'(0)-1$. Let $\mathcal{S}$ denote the subclass of $\mathcal{A}$ which consists of univalent functions in $\ID$ and let
$\mathcal{U}(\lambda),$ $0<\lambda\leq 1,$ denote the set of all  $f\in {\mathcal A}$ satisfying
the condition
\be\label{eq 1}
\left |\left (\frac{z}{f(z)} \right )^{2}f'(z)-1\right | < \lambda \qquad (z\in \ID).
\ee
For $\lambda=1$ we put $\mathcal{U}(1)=\mathcal{U}$. More about these classes can be found  in
\cite{OP-01,OP_2011,OPW_2016,OT_2019,TTV}.

\medskip

In \cite{OPW_2016} it was claimed that all functions $f$ from $\mathcal{U}(\lambda)$ satisfy
\begin{equation}\label{subord}
\frac{f(z)}{z}\prec \frac{1}{(1+z)(1+\lambda z)}.
\end{equation}
Here "$\prec$" denotes the usual subordination, i.e., $F(z)\prec G(z)$, for $f$ and $G$ being analytic functions in $\ID$, means that there exists function $\omega(z)$, also analytic in $\|D$, such that $\omega(0)=0$ and $|\omega(z)|<1$ for all $z\in \ID$. Recently, in \cite{lipon}, the author gave a counterexample that subordination \eqref{subord} is not necessarily satisfied by all functions from $\mathcal{U}(\lambda)$.

\medskip

For the functions $f$ from   $\mathcal{U}(\lambda)$ satisfying subordination \eqref{subord} we have
\be\label{eq 2}
\frac{f(z)}{z}= \frac{1}{(1-\omega(z))(1-\lambda \omega(z))},
\ee
where $\omega$ is a Schwarz function, i.e., it is analytic in $\ID$, $\omega(0)=0$ and $|\omega(z)|<1, z\in \ID$. Let's put
$$\omega(z)=c_{1}z+c_{2}z^{2}+\cdots.$$
Later on we will use the fact due to Schur \cite{schur} that $|c_2|\le1-|c_1|^2$ (can be found also in Carlson's work \cite{carlson}).

\medskip

Further, the inequality \eqref{eq 1} for the function $f$ from $\U(\lambda)$
can be rewritten in the following, equivalent, form
$$ \left|\frac{z}{f(z)}-z \left(\frac{z}{f(z)}\right)'-1 \right|<\lambda \qquad (z\in \ID)$$
and further
$$ \left|\frac{z}{f(z)}-z \left(\frac{z}{f(z)}\right)'-1 \right|\le \lambda|z|^{2}\qquad (z\in \ID).$$
From here, after some calculations we  obtain
$$|(1+\lambda)c_{2}-\lambda c_{1}^{2}+(2(1+\lambda )c_{3}-4\lambda c_{1}c_{2})z+\cdots|\leq\lambda$$
for all $z\in \ID,$ and next,
\be\label{eq 3}
|(1+\lambda)c_{2}-\lambda c_{1}^{2}|\leq\lambda,\quad |2(1+\lambda )c_{3}-4\lambda c_{1}c_{2}|\leq \lambda-\frac{1}{\lambda}|(1+\lambda)c_{2}-\lambda c_{1}^{2}|^{2},
\ee
for all $z\in \ID.$ The last inequality follows from the result of Carlson for the second coefficient of Schwarz functions cited above.

\medskip
If
 $f\in \mathcal{S}$ and
\be\label{eq 4}f(z)=z+a_{2}z^{2}+a_{3}z^{3}+\cdots,
\ee
then the inverse of $f$ has an expansion
\be\label{eq 5}
f^{-1}(w)=w+A_{2}w^{2}+A_{3}w^{3}+\cdots
\ee
near the origin (or precisely at least in $|w|<\frac{1}{4}$).
By using the identity $f(f^{-1})=w$ and the representations for the functions
$f$ and $f^{-1}$, we can obtain the next relations
\be\label{eq 6}
\begin{array}{l}
A_{2}=-a_{2}, \\
A _{3}=-a_{3}+2a_{2}^{2} , \\
A_{4}= -a_{4}+5a_{2}a_{3}-5a_{2}^{3}.
\end{array}
\ee

\medskip

The main result of this paper are the sharp upper bounds for the modulus of these three initial coefficients of $f^{-1}$.

\medskip
\bthm\label{13-th 1}
Let $f\in\mathcal{U}(\lambda)$, $0<\lambda \leq 1$, satisfies subordination \eqref{subord},   and let $f$ and $f^{-1}$ be given by \eqref{eq 4} and \eqref{eq 5}, respectively.
Then
\[
\begin{array}{l}
|A_{2}|\leq 1+\lambda,\\
|A_{3}|\leq1+3\lambda +\lambda^{2},\\
|A_{4}|\leq (1+\lambda)(1+5\lambda +\lambda^{2}).
\end{array}
\]
All these results are best possible.
\ethm

\medskip

\begin{proof}
For $f\in\mathcal{U}(\lambda)$, from the relation \eqref{eq 2} we have (see \cite{obrad-1}, \cite{OPW_2016})
\be\label{eq 8}
\sum_{n=1}^{\infty}a_{n+1}z^{n}=\sum_{n=1}^{\infty}\frac{1-\lambda^{n+1}}{1-\lambda}\omega^{n}(z).
\ee
If we put $\omega(z)=c_{1}z+c_{2}z^{2}+\cdots$, then from \eqref{eq 8} by comparing the coefficients we obtain
\be\label{eq 9}
\left\{\begin{array}{l}
a_{2}=(1+\lambda)c_{1}, \\
a _{3}=(1+\lambda)c_{2}+(1+\lambda+\lambda ^{2})c_{1}^{2} , \\
a_{4}= (1+\lambda)c_{3}+2(1+\lambda+\lambda ^{2})c_{1}c_{2}+(1+\lambda+\lambda ^{2}+\lambda^{3})c_{1}^{3}.
\end{array}
\right.
\ee
Using \eqref{eq 6} and  \eqref{eq 9} we also have
\be\label{eq 10}
\left\{\begin{array}{l}
A_{2}=-(1+\lambda)c_{1}, \\
A _{3}=-(1+\lambda)c_{2}+(1+3\lambda+\lambda ^{2})c_{1}^{2} , \\
A_{4}= -(1+\lambda)c_{3}+(3+8\lambda+3\lambda ^{2})c_{1}c_{2}
-(1+\lambda)(1+5\lambda+\lambda ^{2})c_{1}^{3}.
\end{array}
\right.
\ee
Since $|c_{1}|\leq1$ and $|c_{2}|\leq 1-|c_{1}|^{2}$,  from \eqref{eq 10} we receive
$$|A_{2}|\leq1+\lambda$$
 and
\[
\begin{split}
|A _{3}|&\leq(1+\lambda)|c_{2}|+(1+3\lambda+\lambda ^{2})|c_{1}|^{2}\\
&\leq(1+\lambda)(1-|c_{1}|^{2})+(1+3\lambda+\lambda ^{2})|c_{1}|^{2}\\
&\leq(1+\lambda)+(2\lambda+\lambda ^{2})|c_{1}|^{2}\\
&\leq 1+3\lambda +\lambda^{2} .
\end{split}
\]
Also, from \eqref{eq 10} we obtain
$$A_{4}=-\frac{1}{2}\left[2(1+\lambda )c_{3}-4\lambda c_{1}c_{2}
-6(1+\lambda)c_{1}((1+\lambda)c_{2}  - \lambda c_{1}^{2})  +  2(1+\lambda)^{3}c_{1}^{3}\right],$$
and from here, by applying \eqref{eq 3},
\[
\begin{split}
|A_{4}|&\leq \frac{1}{2}\left[|2(1+\lambda )c_{3}-4\lambda c_{1}c_{2}|
+6(1+\lambda)|c_{1}||(1+\lambda)c_{2}-\lambda c_{1}^{2}|+2(1+\lambda)^{3}|c_{1}|^{3}\right]\\
&\leq \frac{1}{2}\left[\lambda -\frac{1}{\lambda}|(1+\lambda)c_{2}-\lambda c_{1}^{2}|^{2}
+6(1+\lambda)|c_{1}||(1+\lambda)c_{2}-\lambda c_{1}^{2}|+2(1+\lambda)^{3}|c_{1}|^{3}\right]\\
&=\frac{1}{2}\left[\lambda -\frac{1}{\lambda}t^{2}+6(1+\lambda)|c_{1}|t+2(1+\lambda)^{3}|c_{1}|^{3}\right]\\
&=: \frac{1}{2}h(t),
\end{split}
\]
where $t=|(1+\lambda)c_{2}-\lambda c_{1}^{2}|$ and $0\leq t\leq \lambda,$ since
\[ |(1+\lambda)c_{2}-\lambda c_{1}^{2}| \le  (1+\lambda)|c_{2}| + \lambda |c_{1}|^{2} \le  (1+\lambda)(1-|c_1|^2) + \lambda |c_{1}|^{2} = \lambda.\]

\medskip
As for the maximal value of function $h$, we consider two cases:

\medskip
\underline{Case 1:} When $ 0\leq |c_{1}|\leq \frac{1}{3(1+\lambda)}$ the function $h$ attains its maximum for $t_0=3(1+\lambda)\lambda|c_{1}|$ and we have
$$h(t_0)\leq \lambda+27\lambda(1+\lambda)^{2}|c_{1}|^{2}+2(1+\lambda)^{3}|c_{1}|^{3}
\leq 4\lambda +\frac{2}{27},$$
i.e.,
$$|A_{4}|\leq 2\lambda +\frac{1}{27}.$$

\medskip

\underline{Case 2:} For $ \frac{1}{3(1+\lambda)}\leq |c_{1}|\leq 1 $, the function $h$ attains its maximum for $t=\lambda$ and we have
$$h(t)\leq 6(1+\lambda)\lambda|c_{1}|+2(1+\lambda)^{3}|c_{1}|^{3}\leq 2(1+\lambda)(1+5\lambda+\lambda^{2}),$$
when $0\leq t\leq \lambda.$
So,
$$|A_{4}|\leq (1+\lambda)(1+5\lambda+\lambda^{2}).$$

\medskip
From cases 1 and 2, since $(1+\lambda)(1+5\lambda+\lambda^{2}) > 2\lambda +\frac{1}{27}$ when $0<\lambda\le1$, we receive the estimate for $|A_4|$.

\medskip
For the proof of sharpness of the theorem, let consider the function
$$w=f_{\lambda}(z)=\frac{z}{(1-z)(1-\lambda z)}.$$
Then
\be\label{eq 11}
z=f_{\lambda}^{-1}(w)=w-(1+\lambda)w^{2}+ (1+3\lambda +\lambda^{2})w^{3}- (1+\lambda)(1+5\lambda +\lambda^{2})w^{4}-\cdots,
\ee
which shows that our results are the best possible.
\end{proof}

\medskip
Note that for $\lambda=1$ in Theorem 1 we have the estimates for class $\mathcal{U}$ and in that case the inverse of the Koebe function is extremal,
as for the class $\mathcal{S}$ (see, for example Goodman's book, Vol II, p.205, \cite{Go}).

\medskip

In the next theorem we study the Fekete-Szeg\H{o} functional for the inverse functions of the class $\mathcal{U}(\lambda)$. Namely, we have

\medskip

\bthm\label{13-th 1} For the inverse functions of functions from $\mathcal{U}(\lambda )$, $0<\lambda\leq 1$, satisfying subordination \eqref{subord}, we have
$$ |A_{3}-\mu A_{2}^{2}|\leq \lambda +|1-\mu|(1+\lambda)^{2},$$
where $\mu$ is complex number. The result is sharp for $ 0\leq\mu\leq1$.
\ethm

\medskip

\begin{proof}
From the relations \eqref{eq 10} and \eqref{eq 3} we obtain
\[
\begin{split}
|A_{3}-\mu A_{2}^{2}|&=|-(1+\lambda)c_{2}+(1+3\lambda+\lambda^{2})c_{1}^{2}-\mu (1+\lambda)^{2}c_{1}^{2}|\\
&=|-((1+\lambda)c_{2}-\lambda c_{1}^{2})+(1-\mu)(1+\lambda)^{2}c_{1}^{2}|\\
&\leq |(1+\lambda)c_{2}-\lambda c_{1}^{2}|+|1-\mu|(1+\lambda)^{2}|c_{1}|^{2}\\
&\leq \lambda +|1-\mu|(1+\lambda)^{2}.
\end{split}
\]
The sharpness of the estimate in the case when $0\leq\mu\leq1$ follows from the function $f_{\lambda}^{-1}$ defined by \eqref{eq 11}.
\end{proof}

\medskip

\section*{Funding}
The authors did not receive any funding for the research published in this article.

\medskip

\section*{Conflict of Interest}
The authors declare that they have no conflict of interest.

\medskip

\section*{Ethical approval}
This article does not contain any studies with human participants or animals performed by any of the authors.

\end{document}